\documentclass[a4paper,11pt]{article}

\title{Symplectic spreads, planar functions and mutually unbiased bases}
\author{Kanat Abdukhalikov\\
Department of Mathematical Sciences\\
UAE University, PO Box 15551, Al Ain, UAE\\
and \\  
Institute of Mathematics, 
Pushkin Str 125, Almaty 050010, Kazakhstan\\
abdukhalik@uaeu.ac.ae, abdukhalikov@math.kz}

\date{ }

\usepackage{amssymb}

\begin{document}

\maketitle

\newcommand{\Fp}{\mbox{$\mathbb{F}_p$}}          
\newcommand{\Fq}{\mbox{$\mathbb{F}_q$}}          
\newcommand{\Fqt}{\mbox{$\mathbb{F}_q$}}         
\newcommand{\Fqq}{\mbox{$\mathbb{F}_q^*$}}     
\newcommand{\Zf}{\mbox{$\mathbb{Z}_4$}}           
\newcommand{\C}{\mbox{$\mathbb{C}$}}               
\newcommand{\R}{\mbox{$\mathbb{R}$}}               
\newcommand{\T}{\mbox{$\mathcal{T}$}}              

\newtheorem{theorem}{Theorem}[section]
\newtheorem{lemma}[theorem]{Lemma}
\newtheorem{corollary}[theorem]{Corollary}
\newtheorem{proposition}[theorem]{Proposition}

\begin{abstract}
In this paper we give explicit descriptions of complete sets of mutually unbiased bases 
(MUBs) and orthogonal decompositions 
of special Lie algeb\-ras $sl_n(\C)$ obtained from commutative and symplectic semifields, and from 
some other non-semifield symplectic spreads. 
Relations between various constructions are also studied. 
We show that the automorphism group of a complete set of MUBs
is isomorphic to the automorphism group of the corresponding orthogonal 
decomposition of the Lie algebra  $sl_n(\C)$.  
In the case of symplectic spreads  
this automorphism group is determined by the automorphism group of the spread. 
By using the new notion of pseudo-planar functions over fields of characteristic two we  give 
new explicit constructions of complete sets of MUBs. 
\end{abstract}


\section{Introduction}

Mutually unbiased bases (MUBs) were first studied by Schwinger \cite{Schw} in 1960, 
but the notion itself was defined by Wootters and Fields \cite{Woo} almost 30 years later, 
when they also presented examples.
A set of MUBs in the Hilbert space $\C^n$ is defined as  a set of orthonormal 
bases $\{ B_0, B_1, \dots, B_r \}$ of the space such that the square of the absolute value of the 
inner product $|(x,y)|^2$ is equal to $1/n$ for any two vectors $x$, $y$ from distinct bases. 
The notion of MUBs is one of the basic concepts of quantum information theory 
and plays an important role in quantum tomography and state reconstruction \cite{Roy,Woo}. 
It is also valuable for quantum key distribution: the famous Benneth-Brassard secure quantum key 
exchange protocol BB84 and its developments are based on MUBs \cite{Ben}. 
Sequences with low correlations are known to be extremely useful in the design of radar and 
communication systems. 
In 1980 Alltop \cite{All} presented examples of complex sequences with low periodic correlations. 
Later it was noted that these sequences are examples of MUBs. 
On the other hand, in the 1980’s  Kostrikin et al. defined and studied orthogonal decompositions of 
complex Lie algebras. 
The notion of orthogonal decompositions originated in the pioneering work of Thompson, 
who discovered an orthogonal decomposition of the Lie algebra of type $E_8$  and used it for 
construction of a sporadic finite simple group, nowadays called the Thompson group. 
Orthogonal decompositions turn out to be interesting not only for their inner geometric structures, 
but also for their interconnections with other areas of mathematics \cite{Kos}. 
In 2007 Boykin et al. \cite{Boy} discovered a connection between MUBs and 
orthogonal decompositions of Lie algebras. It was found that the existence of a complete set of 
MUBs is equivalent to finding an orthogonal decomposition of the complex Lie algebra $sl_n(\C)$. 
Recently it was discovered that MUBs are very closely related or even equivalent to other problems in 
various parts of mathematics, such as algebraic combinatorics, finite geometry, discrete mathematics, 
coding theory, metric geometry, sequences, and spherical codes. 

There is no general classification of MUBs. The main open problem in this area is to construct a maximal 
number of MUBs for any given $n$. It is known that the maximal set of MUBs of $\C^n$ 
consists of at most $n+1$ bases, and sets attaining this bound are called complete sets of MUBs. 
Constructions of complete sets are known only for prime power dimensions. Even for the smallest non-prime 
power dimension six the problem of finding a maximal set of MUBs is extremely hard  \cite{Kos} 
and remains open after more than 30 years. 

Several kinds of constructions of MUBs are available \cite{All,Ban,Gow,Iva,Kan4,Klap,Woo}. 
As a matter of fact, essentially there are only three types of constructions up to now: 
constructions associated with symplectic spreads, using planar functions over fields of odd characteristic, 
and Gow's construction in \cite{Gow} 
of a unitary matrix (although it seems it is isomorphic to the classical one). 

On the other hand one can consider MUBs in the Euclidean space $\R^n$ (real MUBs), 
and in this case the upper bound for a maximal set is $n/2 + 1$. Constructions of real MUBs are strictly 
connected to algebraic coding theory (optimal Kerdock type codes) and the notion of extremal line-sets 
in Euclidean spaces \cite{Ab4,Cal}. 
LeCompte et al. \cite{LeC} characterized collections of real MUBs in 
terms of association schemes. 
In \cite{Ab4}  connections of real MUBs with binary and quaternary 
Kerdock and Preparata codes \cite{Ham,Nech}, and association schemes were studied.  
We also note that one can come to Kerdock and Preparata codes through integral lattices 
associated with orthogonal decompositions of Lie algebras $sl_n(\C)$ \cite{Ab2,Ab3}. 

In this paper we give explicit descriptions of complete sets of MUBs and orthogonal decompositions 
of special Lie algebras $sl_n(\C)$ obtained from commutative and symplectic semifields, and from 
some other non-semifield symplectic spreads. 
We provide direct formulas to construct MUBs from semifields. 
We show that automorphism groups of complete sets of MUBs and corresponding orthogonal 
decompositions of Lie algebras  $sl_n(\C)$ are isomorphic, and in the case of symplectic spreads  
these automorphism groups are determined by the automorphism groups of those spreads. 
{\rm Aut}omorphism groups are important invariants, therefore they can be used to show inequivalence 
of the various known types of MUBs.  
Planar functions over fields of odd characteristics also lead to constructions of MUBs. 
There are no planar functions over fields of characteristic two. 
However, Zhou \cite{Zhou} proposed a new notion of ``planar" functions over fields of characteristic two. 
Based on this notion,  we propose new constructions of MUBs (but we prefer to call these functions 
pseudo-planar). We also propose a generalization of the notion of pseudo-planar functions 
for arbitrary characteristic. 

MUBs can be constructed using different objects. 
We study their mutual relations. Connections between them are given in the following road map:  

\bigskip

$\put(0,160){\framebox(80,30){\shortstack{Complete\\set of MUBs}}}
\put(225,160){\framebox(80,30){\shortstack{Orth. Decom. of \\  Lie alg.  $sl_{n}(\mathbb{C})$  }}}
\put(0,100){\framebox(80,30){\shortstack{Planar function \\(Pseudo-planar)}}}
\put(230,100){\framebox(80,30){\shortstack{Symplectic\\Spread}}}
\put(0,40){\framebox(80,30){\shortstack{\shortstack{Commutative\\Presemifield}}}}
\put(225,40){\framebox(80,30){\shortstack{Symplectic\\Presemifield}}}
\put(80,175){\vector(1,0){145}} 
\put(225,175){\vector(-1,0){145}}
\put(80,55){\vector(1,0){145}} 
\put(225,55){\vector(-1,0){145}}
\put(40,70){\vector(0,1){30}} 
\put(40,130){\vector(0,1){30}}
\put(265,70){\vector(0,1){30}} 
\put(265,130){\vector(0,1){30}}
$

This map shows that  starting from a commutative presemifield one can construct consecutively 
(pseudo-)planar functions and then MUBs.  
On the other hand, one can start from a commutative presemifield and then construct consecutively 
symplectic presemifield, symplectic spread, orthogonal decomposition of Lie algebra and 
finally MUBs. We will show that in fact our diagram is ``commutative": 
we will not get new MUBs moving from one construction to others. 
Note that there are planar functions not coming from presemifields, 
symplectic spreads not related to semifields \cite{Ball,Kan,Kan8,Kan9,Kan2,Tits} 
and orthogonal decompositions not related to symplectic spreads. 


\section{{\rm Aut}omorphism groups}

The Lie algebra $L=sl_n(\C)$ is the algebra of $n\times n$ traceless matrices over $\C$, 
where the operation of multiplication is given by the commutator of matrices: $[A,B]=AB-BA$. 
A subalgebra $H$ of a Lie algebra $L$ is called a Cartan subalgebra if it is nilpotent and equal to its 
normalizer, which is the set of those elements $X$ in $L$ such that $[X,H] \subseteq H$.  
In case of $L=sl_n(\C)$ Cartan subalgebras are maximal abelian subalgebras and they are 
conjugate under automorphisms of the Lie algebra. In particular, they are all conjugate to 
the standard Cartan subalgebra, consisting of all traceless diagonal matrices.
A decomposition of a simple Lie algebra  $L$ into a direct sum of Cartan subalgebras  
$$L = H_0 \oplus H_1 \oplus \cdots \oplus H_n$$
is called an orthogonal decomposition \cite{Kos}, if the subalgebras $H_i$ are pairwise orthogonal with 
respect to the Killing form $K(A,B)$ on $L$. Recall that the Killing form on a Lie algebra $L$ 
is defined by $K(A,B)={\rm Tr} ({\rm ad} A\cdot  {\rm ad} B )$, where the operator 
${\rm ad} A : L \rightarrow L$ is given by ${\rm ad} A\: (C) = [A,C]$ and ${\rm Tr}$ is the trace. 
The Killing form is symmetric and non-degenerate on $L$.  
In the case of $L=sl_n(\C)$ we have 
$$K(A,B) = 2n{\rm Tr} (AB).$$ 

The adjoint operation $*$ on the set of $n\times n$ complex matrices is given by $A^* = \overline{A}^t$ 
(conjugate transpose). A Cartan subalgebra is called closed under the adjoint operation if $H^* =H$.   

\begin{theorem}[\cite{Boy}, Theorem 5.2]
\label{tiep}
Complete sets of MUBs in $\C^n$ are in one-to-one correspondence with orthogonal 
decompositions of Lie algebra  $sl_n(\C)$ such that all Cartan subalgebras in this decomposition 
are closed under the adjoint operation.  
\end{theorem}

This correspondence is given in the following way. If ${\mathcal B} = \{B_0, B_1, \dots, B_n \}$ 
is a complete set of MUBs then the associated Cartan subalgebra $H_i$ consists of all 
traceless matrices that are diagonal with respect to the basis $B_i$ 
(in other words, vectors of the basis $B_i$ are common eigenvectors of all matrices in $H_i$). 

The automorphism group ${\rm Aut}({\mathcal D})$ \cite{Kos} of an orthogonal decomposition ${\mathcal D}$ 
of a Lie algebra $L$  consists of all automorphisms of $L$ preserving ${\mathcal D}$:
$${\rm Aut}({\mathcal D}) = \{ \varphi \in {\rm Aut}(L) \mid (\forall i) (\exists j) \ \varphi(H_i)=H_j\}.$$

Recall that ${\rm Aut}(L) = {\rm Inn}(L)\cdot \langle T \rangle$, where ${\rm Inn}(L) \cong PSL_n(\C)$ is 
the group of all inner automorphisms of $L$ and $T$ is the outer automorphism 
$A \mapsto -A^t$. 

Let ${\mathcal B} = \{B_0, B_1, \dots, B_n \}$ be a complete set of MUBs of $\C^n$. 
With any orthonormal basis  $B_i$ we can associate an orthoframe $Ort(B_i)$ of 1-spaces generated 
by vectors of $B_i$. We note that if we consider the representations of all bases $B_i$ in 
some fixed standard basis, then the conjugation  map 
$$\tau (x_1, \dots,  x_n) = \overline{x} = (\overline{x_1}, \dots,  \overline{x_n})$$ 
sends one orthoframe $Ort(B_i)$ to other orthoframe. Define 
$$PGL_n(\C)^+ = PGL_n(\C) \langle \tau \rangle ,$$
$${\rm Aut}({\mathcal B}) = \{ \psi \in PGL_n(\C)^+ \mid 
(\forall i) (\exists j) \ \psi(Ort(B_i))=Ort(B_j)\}.$$

\begin{theorem}
\label{autom}
Let ${\mathcal D}$ be an orthogonal decomposition of the Lie algebra  $sl_n(\C)$ and 
${\mathcal B}$ be the corresponding complete set of MUBs. Then 
$${\rm Aut}({\mathcal D}) \cong {\rm Aut}({\mathcal B}).$$
\end{theorem}

{\em Proof}. 
Let $\varphi \in {\rm Aut}({\mathcal D})$. For any automorphism $\varphi$ of the Lie algebra 
$sl_n(\C)$ we have  
$\varphi = \varphi_X$ or $\varphi = \varphi_X T$, where 
$$\varphi_X (A) = XAX^{-1}$$
for some matrix $X\in SL_n(\C)$, and
$$T(A)= - A^t.$$ 
Let $H$ be a Cartan subalgebra from ${\mathcal D}$ and $B=\{e_1, \dots, e_n\}$  
be the corresponding basis, so for any $h\in H$ we have $he_i = \alpha_i(h)e_i$ 
for some linear map $\alpha_i : H \rightarrow \C$. 
Assume that $\varphi = \varphi_X$. 
Then the vectors $f_i =Xe_i$ generate the orthoframe associated with the 
Cartan subalgebra $\varphi(H)$: 
$$XhX^{-1}f_i = XhX^{-1} Xe_i =X\alpha_i(h)e_i = \alpha_i(h) f_i.$$
Therefore, the matrix $X$ generates an element from ${\rm Aut}({\mathcal B})$. 

Let $\varphi = \varphi_X T$. Note that for $h \in H$ one has $\overline{h}^t =h^*  \in H$ 
by Theorem \ref{tiep}.  
Then the vectors $f_i =X\overline{e_i}$ determine a basis corresponding to Cartan subalgebra $\varphi(H)$: 
$$(\varphi_X Th)f_i = -Xh^t X^{-1} f_i = -X\overline{h^*} X^{-1}f_i = 
-X\overline{h^*} X^{-1}X\overline{e_i} = -X\overline{h^* e_i} = $$
$$-X\overline{\alpha_i(h^*)e_i} = -X\overline{\alpha_i(h^*)} \overline{e_i} = 
-\overline{\alpha_i(h^*)}X \overline{e_i} =-\overline{\alpha_i(h^*)} f_i.$$
Therefore, $X\tau$ generates an element  from ${\rm Aut}({\mathcal B})$.

Conversely, let $\psi \in {\rm Aut}({\mathcal B})$. Suppose first that $\psi$ is generated by the matrix $X\in GL_n(\C)$. 
We can assume that $\det X=1$ (otherwise we multiply $X$ by an appropriate scalar matrix). 
Let $B=\{e_i\}$ be a basis from ${\mathcal B}$ and  
let $\psi$ map the orthoframe $Ort(B)$ to another orthoframe $Ort(B')$. 
Denote $f_i  = Xe_i$.  
Let $H'= \varphi_X(H) = XHX^{-1}$ and $h \in H$. Then 
$$\varphi_X(h)f_i = XhX^{-1}Xe_i = Xhe_i = X\alpha_i(h)e_i = \alpha_i(h)Xe_i = \alpha_i(h)f_i .$$ 
Therefore, the vectors $f_i$ are common eigenvectors for $H'$, so $H'$ is the Cartan subalgebra 
associated with the basis $B'$. Hence $\varphi_X \in {\rm Aut}({\mathcal D})$.  

Now assume that $\psi$ is generated by $X\tau$.  
Then the vectors $f_i  = X\overline{e_i}$ from an orthoframe $B'$ are common eigenvectors 
of matrices of the Cartan subalgebra $\varphi_XT(H)$: 
$$(\varphi_X Th)f_i = -Xh^tX^{-1}X\overline{e_i} = -Xh^t\overline{e_i} = 
-X\overline{h^*}\overline{e_i}=-X\overline{h^*e_i}= $$
$$-X\overline{\alpha_i(h^*)e_i} = -X\overline{\alpha_i(h^*)}\overline{e_i} = 
-\overline{\alpha_i(h^*)}X\overline{e_i} = -\overline{\alpha_i(h^*)}f_i ,$$
so  $\varphi_X T \in {\rm Aut}({\mathcal D})$.  $\Box$

\medskip

All known constructions of complete sets of MUBs were obtained with the help of symplectic 
spreads and planar functions, and the construction from \cite{Gow}. 
Now we recall constructions of orthogonal decompositions of 
Lie algebras  $sl_n(\C)$ associated with symplectic spreads \cite{Kos,Kan6}.

Let $F=\Fq$ be a finite field of order $q$. 
Let $V$ be a vector space over $F$ with the usual dot product $u\cdot v$. 
We can consider $W=V\oplus V$ as a vector space over the prime field $\Fp$ and  
define an alternating bilinear form on $W$ by
\begin{eqnarray}
\label{form}
\langle (u,v),(u',v') \rangle = {\rm tr}(u\cdot v' - v\cdot u'),
\end{eqnarray}
where ${\rm tr}$ is a trace function from $\Fq$ to $\Fp$.

Let $n=|V|$ and let $\{e_w\}$ denote the standard basis of $\C^n$, indexed by elements of $V$. 
Let $\varepsilon \in \C$ be a primitive $p$th root of unity. 
For $u \in V$, the generalized Pauli matrices are defined as the following $n\times n$ matrices: 
\begin{eqnarray*}
X(u) & : & e_w \mapsto e_{u+w} \\
Z(v) & : & e_w \mapsto  \varepsilon^{{\rm tr}(v\cdot w)}e_{w}
\end{eqnarray*}

The matrices 
$$D_{u,v} = X(u)Z(v)$$
form a basis of the space of complex square matrices of size $n\times n$. 
Moreover, the matrices  $D_{u,v}$, $(u,v) \ne (0,0)$, generate the Lie algebra $sl_n(\C)$. 
Note that 
$$[D_{u,v},D_{u',v'}] = 
\varepsilon^{{\rm tr}(v\cdot u')}(1 - \varepsilon^{\langle (u,v),(u',v') \rangle})D_{u+u',v+v'},$$
so $[D_{u,v},D_{u',v'}] =0$ if and only if $\langle (u,v),(u',v') \rangle =0$.

Let $V$ be an $r$-dimensional space over $\Fp$ (so $n=p^r$). 
A symplectic spread of the symplectic $2r$-dimensional space $W=V\oplus V$ over $\Fp$ 
is a family of $n+1$ totally isotropic $r$-subspaces of $W$ 
such that every nonzero point of $W$ lies in a unique subspace. 
Thus such $r$-subspaces are maximal totally isotropic subspaces. 
Let $\Sigma = \{ W_0, W_1, \dots, W_n \}$ be a symplectic spread. Then 
$$L = H_0 \oplus H_1 \oplus \cdots \oplus H_n,$$
$$H_i = \langle D_{u,v} \mid (u,v)\in W_i \rangle ,$$
gives \cite{Kos,Kan6} the corresponding orthogonal decomposition of the Lie algebra $sl_n(\C)$. 

We define 
$$Sp^{\pm}(W) = \{ \varphi \in GL(W) \mid 
(\exists s=\pm 1)(\forall u, v\in W) \ \langle \varphi(u),\varphi(v) \rangle = s\langle u,v \rangle \}$$
$${\rm Aut}(\Sigma) = \{ \varphi \in Sp^{\pm}(W) \mid (\forall i )(\exists j) \ \varphi(W_i)=W_j \}$$

\begin{theorem}[\cite{Kos}, Proposition 1.3.3.]
\label{aut-spread}
Let ${\mathcal D}$ be the orthogonal decomposition of the Lie algebra  $sl_n(\C)$ 
corresponding to a symplectic spread $\Sigma$. Then  
$${\rm Aut}({\mathcal D}) = K . {\rm Aut}(\Sigma),$$
where $K$ is the set of all automorphisms  that fix every line $\langle D_{u,v} \rangle$. 
\end{theorem}

In the case of odd characteristic, $K$ is isomorphic to $W=V^2$ and an embedding of $K$ in $PSL_n(\C)$ 
is given through conjugations by matrices $D_{u,v}$: 
$$D_{u,v}D_{a,b}D_{u,v}^{-1} = \varepsilon^{-\langle (u,v), (a,b) \rangle}D_{a,b}.$$

The outer automorphism $T$ is given by:
$$T(D_{a,b}) = -\varepsilon^{-{\rm tr}(a\cdot b)} D_{-a,b},$$
therefore in the case of even characteristic, we have $K \cong V^2 . \ \mathbb{Z}_2$. 

Any symplectic spread $\Sigma$ determines a translation plane $A(\Sigma)$. 
The collineation group of  $A(\Sigma)$ is a semidirect product of the translation group $V^2$ and 
the tran\-slation complement. 
However, the group extension in Theorem \ref{aut-spread} can be nonsplit (see Section \ref{des2}).


\section{Symplectic spreads in odd characteristics and MUBs}
\label{odd-char}
In this section we show how to construct a complete set of MUBs directly from symplectic spreads and semifields. 
Throughout this section $F=\Fq$ denotes a finite field of odd order $q=p^r$. 
Let $\omega \in \C$ be a primitive $p$th root of unity. 

\begin{lemma}
\label{lemma1}
Let $\Sigma$ be a symplectic spread of $W=V \oplus V$, and let $h: V \rightarrow V$ be a mapping 
such that $\{(u,h(u)) \mid u \in V\}$ is a maximal totally isotropic subspace. 
Then ${\rm tr}(u\cdot h(w)) = {\rm tr} (h(u)\cdot w)$ for all $u\in V$, $w\in V$. 
\end{lemma}

{\em Proof}. ${\rm tr}(u\cdot h(w) - h(u)\cdot w) = \langle ( (u,h(u)),(w,h(w)) \rangle =0$. $\Box$

\begin{lemma}
\label{lemma2}
Let $\Sigma$ be a symplectic spread of $W=V \oplus V$, and let $h: V \rightarrow V$ be a linear mapping 
such that $\{(u,h(u)) \mid u \in V\}$ is a maximal totally isotropic subspace. 
Then for any $v\in V$, the vector 
$$b_{h,v} =  \sum_{w\in V} \omega^{{\rm tr}(\frac{1}{2}w\cdot h(w) + v\cdot w)} e_w$$
is an eigenvector of $D_{u,h(u)}$ for all $u\in V$. 
\end{lemma}

{\em Proof}. Indeed, 

\begin{eqnarray*}
D_{u,h(u)}(b_{h,v}) 
& = &  \sum_{w\in V} \omega^{{\rm tr}(\frac{1}{2}w\cdot h(w) + v\cdot w + h(u)\cdot w)} e_{w+u} \\
& = &  \sum_{w\in V} 
           \omega^{{\rm tr}(\frac{1}{2}(w-u)\cdot h(w-u) + v\cdot (w-u) + h(u)\cdot (w-u))} e_{w} \\
& = &  \sum_{w\in V} \omega^{{\rm tr}(\frac{1}{2}w\cdot h(w) + v\cdot w) + 
                         \frac{1}{2}{\rm tr}(h(u)\cdot w - u\cdot h(w)) -  {\rm tr}(\frac{1}{2} u\cdot h(u) + v\cdot u)} e_{w} \\
& = & \omega^{ - {\rm tr}(\frac{1}{2} u\cdot h(u) + v\cdot u)} b_{h,v} 
\end{eqnarray*}
by Lemma \ref{lemma1}. $\Box$

Lemma \ref{lemma2} allows us to construct the complete set of MUBs 
corresponding to the orthogonal decomposition of the Lie algebra $sl_n(\C)$ 
obtained from a symplectic spread.  


Symplectic spreads can be constructed using semifields. A finite presemifield is a ring with no 
zero-divisors, and with left and right distributivity \cite{Dem}. 
A presemifield with multiplicative identity is called a semifield.  
A finite presemifield can be obtained from a finite field $(F,+,\cdot)$ by introducing a new 
product operation $*$, so it is denoted by $(F,+,*)$. 
Every presemifield determines a spread $\Sigma$ consisting of subspaces 
$(0,F)$ and $\{ (x,x*y) \mid x\in F\}$, $y\in F$.  A presemifield is called symplectic if the 
corresponding spread is symplectic with respect to some alternating form \cite{Kan3} 
(so the spread might not be symplectic with respect to other forms). 

Two presemifields $(F,+,*)$ and $(F,+,\star)$ 
are called isotopic if there exist three bijective linear mappings 
$L$, $M$, $N: F \rightarrow F$ such that 
$$L(x*y) = M(x) \star N(y)$$
for any $x, y \in F$. \  If $M=N$ then the presemifields are called strongly isotopic. 
Every presemifield is isotopic to a semifield. 
Isotopic semifields determine isomorphic planes. 

Lemma \ref{lemma2} implies

\begin{theorem}
\label{symp}
Let $(F,+,\circ )$ be a finite symplectic presemifield of odd characteristic. Then the following set 
forms a complete set of MUBs: 
$$B_{\infty}= \{ e_w \mid w \in F \}, \quad B_{m} = \{ b_{m,v} \mid v \in F \}, \quad m\in F ,$$ 
$$b_{m,v} = \frac{1}{\sqrt{q}} \sum_{w \in F} \omega^{{\rm tr}(\frac{1}{2}w \cdot (w\circ m) + v\cdot w)}e_w.$$
\end{theorem}

A function $f : F \rightarrow F$ is called planar if 
\begin{eqnarray}
\label{planarf}
x \mapsto f(x+a)-f(x) ,
\end{eqnarray}
is a permutation of $F$ for each $a\in F^*$. 
Any planar function over $F$  allows us to construct a complete set of MUBs.  

\begin{theorem}[\cite{Ding,Kan4,Roy}]
\label{planar}
Let $F$ be a finite field of odd order $q$  and $f$ be a planar function. 
Then the following  forms a complete set of MUBs:
 $$B_{\infty}= \{ e_w \mid w \in F \}, \quad B_{m} = \{ b_{m,v} \mid v \in F \}, \quad m\in F ,$$ 
$$b_{m,v} = \frac{1}{\sqrt{q}} \sum_{w \in F} \omega^{{\rm tr}(\frac{1}{2}mf(w) + vw)}e_w.$$
\end{theorem}

Every commutative presemifield of odd order corresponds to a quadratic planar 
polynomial, and vice versa. If $(F,+,*)$ is a commutative presemifield then $f(x)=x*x$ 
is a planar function, and one can use Theorem \ref{planar} for constructing MUBs. 
Up to now, there is only one known type of planar function \cite{Coul}, 
which is not related to semifields: 
$$f(x) = x^{(3^k+1)/2},$$
where $q=3^r$,  $\gcd(k,2r)=1$, $k \not\equiv \pm 1 \!\!\!\pmod{2r}$. 

On the other hand, from a commutative presemifield one can construct a symplectic presemifield 
(using Knuth's cubical array method \cite{Knu}) and then construct  consecutively 
a symplectic spread, an orthogonal decomposition of the Lie algebra $sl_n(\mathbb{C})$ and 
finally MUBs. Below we show that we will not get  new MUBs by these operations 
(a similar statement is true for even characteristic, see Section \ref{even-char}).

Let $f(x)= \sum_{i\leq j}a_{ij}x^{p^i +p^j}$ be a quadratic planar polynomial.  
Then the corresponding commutative presemifield $(F, +, * )$ is given by:
$$x * y = \frac{1}{2}(f(x+y)-f(x)-f(y)) = 
\frac{1}{2} \sum_{i\leq j}a_{ij}(x^{p^i}y^{p^j} + x^{p^j}y^{p^i}).$$ 
It defines a spread $\Sigma$ consisting of subspaces 
$(0,F)$ and $\{ (x,x*y) \mid x\in F\}$, $y\in F$.   
Starting from $\Sigma$ we will construct a symplectic spread using the Knuth cubical arrays method 
(\cite{Kan3}, Proposition 3.8). 
The dual spread $\Sigma^d$ is a spread of the dual space of $F\oplus F$. 
We can identify that dual space with $F\oplus F$ by using the alternating form (\ref{form}). 
For each $y\in F$ we have to find all $(u,v)$ such that $\langle (x,x*y), (u,v) \rangle = 0$ 
for all $x \in F$. We have 
\begin{eqnarray*}
\langle (x,x*y), (u,v) \rangle 
& = & {\rm tr} (xv - u (x*y)) \\
& = & {\rm tr} (xv - u \frac{1}{2} \sum_{i\leq j}a_{ij}(x^{p^i}y^{p^j} + x^{p^j}y^{p^i}) ) \\
& = & {\rm tr} (xv - \frac{1}{2} \sum_{i\leq j}a_{ij}^{p^{r-i}}u^{p^{r-i}}xy^{p^{j-i}} -  
                    \frac{1}{2} \sum_{i\leq j}a_{ij}^{p^{r-j}}u^{p^{r-j}}xy^{p^{r+i-j}}  ) \\
& = & {\rm tr} (x (v - \frac{1}{2} \sum_{i\leq j}a_{ij}^{p^{r-i}}u^{p^{r-i}}y^{p^{j-i}} -  
                    \frac{1}{2} \sum_{i\leq j}a_{ij}^{p^{r-j}}u^{p^{r-j}}y^{p^{r+i-j}}  )) ,
\end{eqnarray*}
which implies 
$$v = \frac{1}{2} \sum_{i\leq j}a_{ij}^{p^{r-i}}u^{p^{r-i}}y^{p^{j-i}} +  
                    \frac{1}{2} \sum_{i\leq j}a_{ij}^{p^{r-j}}u^{p^{r-j}}y^{p^{r+i-j}} .$$ 
Therefore, $\Sigma^d$ corresponds to the presemifield $(F, +, \bullet )$ defined by
$$u \bullet y =  \frac{1}{2} \sum_{i\leq j}a_{ij}^{p^{r-i}}u^{p^{r-i}}y^{p^{j-i}} +  
                     \frac{1}{2} \sum_{i\leq j}a_{ij}^{p^{r-j}}u^{p^{r-j}}y^{p^{r+i-j}} .$$ 
Then the presemifield $(F, +, \circ )$ with multiplication 
$$x \circ y = y \bullet x = \frac{1}{2} \sum_{i\leq j}a_{ij}^{p^{r-i}} x^{p^{j-i}} y^{p^{r-i}}  +  
                     \frac{1}{2} \sum_{i\leq j}a_{ij}^{p^{r-j}} x^{p^{r+i-j}} y^{p^{r-j}} $$
defines a symplectic spread $\Sigma^{d*}$.  
It is straightforward to check directly that the spread $\Sigma^{d*}$ with subspaces 
$(0,F)$ and $\{ (x,x\circ y) \mid x\in F\}$, $y\in F$, is symplectic with respect to the form  (\ref{form}). 

Then by Theorem \ref{symp} the following bases will form a complete set of MUBs:
$$B_{\infty}= \{ e_w \mid w \in F \}, \quad B_{m} = \{ b_{m,v} \mid v \in F \}, \quad m\in F ,$$ 
$$b_{m,v} = \frac{1}{\sqrt{q}} \sum_{w \in F} \omega^{{\rm tr}(\frac{1}{2}m(w*w) + vw)}e_w,$$
since 
\begin{eqnarray*}
{\rm tr} \left( \frac{1}{2}w \cdot (w \circ m) \right) 
& = &  {\rm tr} \left(\frac{1}{4} \sum_{i\leq j} w a_{ij}^{p^{r-i}} w^{p^{j-i}} m^{p^{r-i}}  +  
                     \frac{1}{4} \sum_{i\leq j} w a_{ij}^{p^{r-j}} w^{p^{r+i-j}} m^{p^{r-j}} \right)  \\
& = &  {\rm tr} \left(\frac{1}{4} \sum_{i\leq j} w^{p^i} a_{ij} w^{p^{j}} m  +  
                     \frac{1}{4} \sum_{i\leq j} w^{p^j} a_{ij} w^{p^{i}} m \right)  \\
& = & {\rm tr} \left( \frac{1}{2}m(w*w) \right).
\end{eqnarray*}

{\bf Remark 1}. 
It is easy to see that strongly isotopic commutative semifields provide equivalent MUBs 
(equivalence under the action of the group $PGL_n(\mathbb{C})^+$). 
Theorems \ref{autom} and \ref{aut-spread} can be used to show inequivalence of MUBs 
obtained from nonisotopic semifields.  Practically they distinguish all known types of MUBs.


\subsection{Desarguesian spreads}

Desarguesian spreads are constructed with the help of finite fields, and the corresponding commutative 
and symplectic semifields are the same: $x*y = x\circ y =xy$, and the planar function is $f(x)=x^2$. 

Godsil and Roy \cite{God} showed that the constructions of Alltop \cite{All}, Ivanovi\'{c} \cite{Iva}, 
Wooters and Fields \cite{Woo}, Klappenecker and R\"{o}tteler \cite{Klap},  and 
Bandyopadhyay, Boykin, Roychowdhury and Vatan \cite{Ban} are all equivalent to 
particular cases of this construction. It seems that Gow's construction is equivalent to this case 
as well.

Theorem \ref{autom}  and \cite{Kos} imply that for the automorphism group of the 
correspon\-ding complete set of MUBs  we have 
$${\rm Aut}(\mathcal{B}) \cong F^2 . (SL_2(q) . \mathbb{Z}_r) . \mathbb{Z}_2 ,$$
where extensions are split.


\subsection{Spreads from Albert's generalized twisted fields}

Let $\rho$ be a nontrivial automorphism of a finite field $F$ such that $-1 \not\in F^{\rho -1}$. 
It means that $F$ has odd degree over the fixed field $F_{\rho}$ of $\rho$. 
Let $V=F$. Then the BKLA \cite{Bad} symplectic spread $\Sigma$ of $W=V\oplus V$ 
consists of the subspace $\{(0,y) \mid y \in F \}$ and subspaces 
$\{(x,mx^{\rho^{-1}} + m^{\rho}x^{\rho}) \mid x \in F \}$, $m\in F$. 
This spread arises from a presemifield given by the operation 
$x\circ m =  mx^{\rho^{-1}} + m^{\rho}x^{\rho}$. 


If we start from the planar function $f(x)=x^{{p^k}+1} = x^{\rho +1}$, then the commutative 
presemifield and  corresponding symplectic presemifield are given by  
$$x*y = \frac{1}{2}(x^{\rho}y + xy^{\rho}),$$
$$x \circ y = \frac{1}{2}(x^{\rho}y + x^{\rho^{-1}}y^{\rho^{-1}}).$$ 




\subsection{Ball-Bamberg-Laurauw-Penttila symplectic spread}

This non-semifield spread \cite{Ball} is obtained from the previous one with the help of a technique 
known as ``net replacement". 
Let $\Sigma$ be the BKLA spread. 
One can change this spread and get a new spread $\Sigma'$ in the following way. 
For all $s \in F^*$ the map $\sigma_s(v,w)=(sv,s^{-1}w)$ is an isometry of $V\oplus V$ 
with respect to the form $\langle \ , \rangle$. 
Denote the group of all $\sigma_s$  by $G$. 
Let $\tau$ be an involution of $V\oplus V$ which switches coordinates:  
$(v,w) \mapsto (w,v)$. Then $\Sigma$ is $G$-invariant. 
For $W_1= \{ (x,x^{\rho} + x^{\rho^{-1}}) \mid x \in F \}$ we consider the $G$-orbit 
$${\mathcal N} = \{ \sigma_s(W_1) \mid s \in F^* \}.$$
Note that ${\mathcal N}$ and $\tau({\mathcal N})$ are $G$-invariant. Now we define 
$${\mathcal N}' = \tau({\mathcal N}).$$
Then 
$$\Sigma' = ( \Sigma \setminus {\mathcal N} ) \cup {\mathcal N}' $$ 
is a symplectic spread  with respect to the form (\ref{form}). 

\begin{lemma}
\label{orbit}
Under the action of the group $G$  the spread $\Sigma$ has the following four orbits: 
$\{(x,0) \mid x \in F \}$; $\{(0,y) \mid y \in F \}$; 
subspaces $\{(x,mx^{\rho^{-1}} + m^{\rho}x^{\rho}) \mid x \in F \}$ with nonsquare elements 
$m\in F^*$; 
and subspaces $\{(x, mx^{\rho^{-1}} + m^{\rho}x^{\rho}) \mid x \in F \}$ with square elements 
$m\in F^*$.
\end{lemma}

{\em Proof}. We have 
\begin{eqnarray*}
\sigma_s(x,mx^{\rho^{-1}} + m^{\rho}x^{\rho})  
& = & (sx,s^{-1}mx^{\rho^{-1}} + s^{-1}m^{\rho}x^{\rho})   \\
& = & (sx,s^{-1-\rho^{-1}}m(sx)^{\rho^{-1}} + s^{-1-\rho}m^{\rho}(sx)^{\rho})  \\
& = & (u,s^{-1-\rho^{-1}}mu^{\rho^{-1}} + s^{-1-\rho}m^{\rho}u^{\rho}) ,
\end{eqnarray*}
where $u=sx$. Therefore, the isometry $\sigma_s$ sends the subspace indexed by $m$  
to the subspace indexed by  $s^{-1-\rho^{-1}}m$. 
Since $s^{-1-\rho^{-1}} = (s^{\rho +1})^{-\rho^{-1}}$, it remains for us to show that 
the set $\{ s^{\rho +1} \mid s \in F^*\}$ is the set of all squares in $F^*$. 
Indeed, assume that the fixed field $F_{\rho}$ of $\rho$ has $q_1$ elements. 
Then $\rho =q_1^k$, $|F|=q_1^t$, $\gcd(t,k)=1$, and $t$ is an odd integer. We have 
$\gcd(q_1^k +1, q_1^t -1) =  2$, so $\{ s^{\rho +1} \mid s \in F^*\}$ is the set of all squares in $F^*$. 
$\Box$

\begin{lemma}
\label{inverse}
Suppose the function $\alpha : F \rightarrow F$ is given by $\alpha(u)=u^{\rho^{-1}} +u^{\rho}$. 
Let $\beta =\alpha^{-1}$ and let the order of $\rho$ be $t$. Then 
$$\beta(v) = \frac{1}{2}(a_0 v + a_1v^{\rho} + \cdots + a_{t-1}v^{\rho^{t-1}}),$$
where $a_{4i} = a_{4i+1} =1$, $a_{4i+2} = a_{4i+3} =-1$, $i \ge 0$, for $t \equiv 1 \pmod{4}$,  
and $a_{4i} = a_{4i+3} =-1$, $a_{4i+1} = a_{4i+2} =1$, $i \ge 0$, for $t \equiv 3 \pmod{4}$.
\end{lemma}

{\em Proof}. We have 
\begin{eqnarray*}
\beta(v) + \beta(v)^{\rho^2}   
& = & \frac{1}{2}(a_0 v + \cdots + a_{t-1}v^{\rho^{t-1}})  \\
&    &  + \frac{1}{2}(a_0v^{\rho^2} + \cdots + a_{t-3}v^{\rho^{t-1}} + a_{t-2}v + a_{t-1}v^{\rho})  \\
& = & \frac{1}{2}((a_0 + a_{t-2})v + (a_1 + a_{t-1})v^{\rho} + 
                            \sum_{k=2}^{t-1} (a_k + a_{k-2})v^{\rho^k}) = v^{\rho}. 
\end{eqnarray*}
Therefore $\beta(v)^{\rho^{-1}} + \beta(v)^{\rho} = v$. 
$\Box$

We use this function $\beta$ for the following 

\begin{corollary}
\label{cor1}
The spread $\Sigma'$ consists of the following subspaces: 
$\{(x,0) \mid x \in F \}$; $\{(0,y) \mid y \in F \}$; 
subspaces $\{(x,mx^{\rho^{-1}} + m^{\rho}x^{\rho}) \mid x \in F \}$ with nonsquare elements $m\in F$; 
and subspaces $\{( x, s\beta(xs)) \mid x \in F \}$ with $s\in F^*/\langle -1\rangle$.
\end{corollary}

{\em Proof}. The collection ${\mathcal N}'$ consists of elements of the form 
$(s^{-1}(u^{\rho^{-1}} + u^{\rho}), su)$. 
We denote $s^{-1}(u^{\rho^{-1}} + u^{\rho}) = x$. 
Then $u^{\rho^{-1}} + u^{\rho} = xs$, $u=\beta(xs)$  and $su=s\beta(xs)$. 
Therefore the corresponding spread subspaces are $\{ (x, s\beta(xs)) \mid x \in F \}$, $s\in F^*$. 
Finally we note that $s$ and $-s$ determine the same subspaces.   
$\Box$

From this corollary we get the following  complete set of MUBs:

$$B_{\infty}= \{ e_w \mid w \in F \}, \ B_0 = \{ \frac{1}{\sqrt{q}} \sum_{w \in F} \omega^{{\rm tr}(vw)}e_w \mid v \in F \},$$
$$B_{m} = \{ \frac{1}{\sqrt{q}} \sum_{w \in F} \omega^{{\rm tr}(m^{\rho}w^{\rho +1} + vw)}e_w 
\mid v \in F \}, \quad m\in F^*, \ m \ {\rm is \ nonsquare},$$ 
$$B_{s} = \{ \frac{1}{\sqrt{q}} \sum_{w \in F} \omega^{{\rm tr}(\frac{1}{2}ws\beta(ws) + vw)}e_w 
\mid v \in F \}, \quad s\in F^*/\langle -1\rangle.$$

We note that the Ball-Bamberg-Laurauw-Penttila symplectic spread is not a semifield spread. 
If $F$ is a field of order 27, then the corresponding plane is the Hering plane and its translation 
complement is isomorphic to $SL_2(13)$. If $|F|=q_1^3$ and $|F_{\rho}|=q_1$ then the corresponding 
plane is the Suetake \cite{Sue} plane and the translation component has order $3(q_1-1)(q_1^3-1)$. 

\subsection{Other known semifields}

In this section we consider examples of known commutative and symplectic semifields. 
They are presented in pairs, so that symplectic semifields produce symplectic spreads 
with respect to the alternating form (\ref{form}).

The Dickson commutative semifield \cite{Dem,Dick} and corresponding Knuth \cite{Knu} 
symplectic presemifields are defined by   
$$(a,b) * (c,d) = (ac + jb^{\sigma}d^{\sigma}, ad + bc),$$
$$(a,b) \circ (c,d) = (ac+bd, ad + j^{\sigma^{-1}}bc^{\sigma^{-1}}),$$
where  $q$ is odd, $j$ be a nonsquare element of $F=\Fq$,   $1\neq \sigma \in {\rm Aut}(F)$ and  
 $V=F \oplus F$. 


The Cohen-Ganley commutative semifield  \cite{Coh} and the corresponding 
symplectic Thas-Payne presemifield \cite{Thas} are defined by 
$$(a,b) * (c,d) = (ac + jbd + j^3(bd)^9, ad + bc + j(bd)^3),$$
$$(a,b) \circ (c,d) = (ac +bd, ad + jbc + j^{1/3}bc^{1/9} + j^{1/3}bd^{1/3}),$$
where  $q=3^r \ge 9$ and  $j \in F$ is nonsquare. 


Ganley commutative semifields \cite{Gan} and corresponding symplectic presemifields 
are defined by 
$$(a,b) * (c,d) = (ac - b^9d -bd^9, ad + bc +b^3d^3),$$
$$(a,b) \circ (c,d) = (ac + bd, ad + bd^{1/3} - b^{1/9}c^{1/9} - b^9c ),$$
where $F=\Fq$, $q=3^t$, and $t \ge 3$ odd . 


The Penttila-Williams sporadic symplectic semifield \cite{Pen} 
$(\mathbb{F}_{3^5} \oplus \mathbb{F}_{3^5}, +,\circ)$ of order $3^{10}$  
and the corresponding commutative semifield are given by  
$$(a,b) \circ (c,d) = (ac + bd, ad + bd^9 +bc^{27}),$$
$$(a,b) * (c,d) = (ac + (bd)^9, ad + bc +(bd)^{27}).$$


During last several years some families of quadratic planar functions were discovered 
\cite{Bier,Buda,Coul,Coul2,Coul3,Ding1,Weng,Zha,Zha2,Zhou2}.


\section{Symplectic spreads in even characteristics and MUBs}
\label{even-char}
In this section we use symplectic spreads in even characteristic for constructing MUBs.
Constructions of MUBs are also given in \cite{Kan4}, where related objects are considered 
over prime fields. We will use Galois rings, so descriptions can be given explicitly 
(especially in the case of presemifields).  
Throughout this section $F=\Fq$ denotes a finite field of even order $q=2^r$. 
Let $\omega \in \C$ be a primitive $4$th root of unity. 

To write corresponding MUBs, we need the Galois ring $R=GR(4^r)$ of cha\-racteristic 4 and cardinality $4^r$. 
We recall some facts about $R$ \cite{Ham,Nech}. 
We have $R/2R \cong \Fq$,  the unit group $R^* = R\setminus 2R$ 
contains a cyclic subgroup $C$ of size $2^r-1$  isomorphic to $\Fqq$. 
The set $\T= \{0\} \cup C$ is called the Teichm\"{u}ller set in $R$. 
Every element $x \in R$ can be written uniquely in the form $x=a+2b$ for $a$, $b \in \T$. 
Then 
$${\rm Tr} (x) = (a+a^2 + \cdots + a^{2^{r-1}}) +   2(b+b^2 + \cdots +b^{2^{r-1}}).$$
Since $R/2R \cong \Fq$, for every element $u \in \Fq$ there exists a corresponding unique 
element $\hat{u} \in \T$, called the Teichm\"{u}ller lift of $u$. 
If $x$, $y$, $z \in \T$ and $z\equiv x+y \!\!\!\pmod{2R}$ then 
$$z = x+y+2\sqrt{xy}.$$ 
The last equation can be written in the form 
$$\hat{w} = \hat{u} + \hat{v} + 2\sqrt{\hat{u}\hat{v}}$$
for elements $w$, $u$, $v \in \Fq$, $w=u+v$. 

Let $(F, +, \circ)$ be a presemifield with multiplication 
$$x\circ y = \sum a_{ij}x^{2^i}y^{2^j}.$$ 
Then we extend this multiplication to $\T \times \T$ in the following way:  
$$\hat{x}\circ \hat{y} =  \sum \widehat{a_{ij}}\hat{x}^{2^i}\hat{y}^{2^j}$$

\begin{lemma}
\label{lemma-tr2}
Let $(F, +, \circ)$ be a symplectic presemifield. Then 

1. ${\rm tr}(x\cdot (z \circ y)) = {\rm tr} (z \cdot (x\circ y))$ for all $x$, $y$, $z\in F$. 

2.  ${\rm Tr} (\hat{x}\cdot (\hat{z} \circ \hat{y})) = 
{\rm Tr}  (\hat{z} \cdot (\hat{x}\circ \hat{y}))$ for all $x$, $y$, $z\in F$. 
\end{lemma}

{\em Proof}. 
1. ${\rm tr}(x\cdot (z \circ y) - z \cdot (x\circ y)) = \langle ( (x,x\circ y),(z,z\circ y) \rangle =0$. 

2. We have 
\begin{eqnarray*}
{\rm tr}(z \cdot (x\circ y) - x\cdot (z \circ y))
& = & {\rm tr} (z \sum a_{ij}x^{2^i}y^{2^j} - x \sum a_{ij}z^{2^i}y^{2^j} ) \\
& = & {\rm tr} (z \sum a_{ij}x^{2^i}y^{2^j} -  \sum a_{ij}^{2^{-i}}x^{2^{-i}}zy^{2^{j-i}} ) \\
& = & {\rm tr} (z (\sum a_{ij}x^{2^i}y^{2^j} -  \sum a_{ij}^{2^{-i}}x^{2^{-i}}y^{2^{j-i}} )) =0, 
\end{eqnarray*}
where  powers of 2 and indices are considered modulo $r$ for convenience. Therefore, 
$$\sum a_{ij}x^{2^i}y^{2^j} =  \sum a_{ij}^{2^{-i}}x^{2^{-i}}y^{2^{j-i}} = 
\sum a_{r-i,j}^{2^{i}}x^{2^{i}}y^{2^{j+i}} = 
\sum a_{r-i,j-i}^{2^{i}}x^{2^{i}}y^{2^{j}}.$$

Hence $a_{ij} = a_{r-i,j-i}^{2^{i}}$ and $\widehat{a_{ij}} = \widehat{a_{r-i,j-i}}^{2^{i}}$. Then 
\begin{eqnarray*}
{\rm Tr} (\hat{z} \cdot (\hat{x}\circ \hat{y}) - \hat{x}\cdot (\hat{z} \circ \hat{y}))
& = & {\rm Tr}  (\hat{z} \sum \widehat{a_{ij}}\hat{x}^{2^i}\hat{y}^{2^j} - 
                          \hat{x} \sum \widehat{a_{ij}}\hat{z}^{2^i}\hat{y}^{2^j} ) \\
& = & {\rm Tr}  (\hat{z} (\sum \widehat{a_{ij}}\hat{x}^{2^i}\hat{y}^{2^j} -  
                          \sum \widehat{a_{ij}}^{2^{-i}}\hat{x}^{2^{-i}}\hat{y}^{2^{j-i}} )) \\
& = & {\rm Tr}  (\hat{z} (\sum \widehat{a_{ij}}\hat{x}^{2^i}\hat{y}^{2^j} -  
                          \sum \widehat{a_{r-i,j}}^{2^{i}}\hat{x}^{2^{i}}\hat{y}^{2^{j+i}} )) \\
& = & {\rm Tr}  (\hat{z} (\sum \widehat{a_{ij}}\hat{x}^{2^i}\hat{y}^{2^j} -  
                          \sum \widehat{a_{r-i,j-i}}^{2^{i}}\hat{x}^{2^{i}}\hat{y}^{2^{j}} )) =0. \ \Box
\end{eqnarray*}

\begin{theorem}
\label{symp2}
Let $(F, +, \circ)$ be a symplectic presemifield. 
Then the following forms a complete set of MUBs:
$$B_{\infty}= \{ e_w \mid w \in F \}, \quad B_{m} = \{ b_{m,v} \mid v \in F \}, \quad m\in F ,$$ 
$$b_{m,v} = \frac{1}{\sqrt{q}} \sum_{w \in F} 
     \omega^{{\rm Tr} (\hat{w}\cdot (\hat{w}\circ \hat{m}) + 2\hat{w}\hat{v})}e_w .$$
\end{theorem}

{\em Proof}. We show that for any $m$, $v\in F$, the vector 
$$d_{m,v} = \sum_{w \in F} \omega^{{\rm Tr} (\hat{w}\cdot (\hat{w}\circ \hat{m}) + 2\hat{w}\hat{v})}e_w$$
is an eigenvector of $D_{u,u\circ m}$ for all $u\in V$. Indeed, 
\begin{eqnarray*}
D_{u,u\circ m}(d_{m,v}) 
& = &  \sum_{w\in F} \omega^{{\rm Tr} (\hat{w}\cdot (\hat{w}\circ \hat{m}) + 
              2\hat{w}\hat{v} + 2\hat{w} \cdot \widehat{u\circ m})}e_{w+u} \\
& = &  \sum_{w\in F} \omega^{{\rm Tr} (\widehat{(w+u)}\cdot (\widehat{(w+u)}\circ \hat{m}) + 
              2(\hat{w} + \hat{u})\hat{v} + 2(\hat{w}+\hat{u}) \cdot (\hat{u}\circ \hat{m}))}e_{w} \\
& = &  \sum_{w\in F} 
            \omega^{{\rm Tr} ((\hat{w}+\hat{u}+ 2\sqrt{\hat{w}\hat{u}}) \cdot 
                            ((\hat{w}+\hat{u}+ 2\sqrt{\hat{w}\hat{u}})\circ \hat{m}) + 
              2(\hat{w} + \hat{u})\hat{v} + 2(\hat{w}+\hat{u})  (\hat{u}\circ \hat{m}))}e_{w} \\
& = &  \sum_{w \in F} 
            \omega^{{\rm Tr} (\hat{w}\cdot (\hat{w}\circ \hat{m}) + 2\hat{w}\hat{v} 
                       + S + \hat{u} (\hat{u}\circ \hat{m}) + 2\hat{u}\hat{v} + 
                                     2\hat{u}(\hat{u}\circ\hat{m}) )}e_w \\
& = &  \omega^{{\rm Tr} (2\hat{u}\hat{v} - \hat{u}(\hat{u}\circ\hat{m}) )} d_{m,v} 
\end{eqnarray*}
where 
\begin{eqnarray*}
{\rm Tr} (S)
& = & {\rm Tr} ((\hat{w}\cdot (\hat{u}\circ\hat{m})  + \hat{u}\cdot (\hat{w}\circ\hat{m}) + 
                                 2\hat{w}\cdot (\hat{u}\circ\hat{m}) )  \\
&  &       +  (2\hat{w}\cdot (\sqrt{\hat{w}\hat{u}} \circ \hat{m}) + 
                     2\sqrt{\hat{w}\hat{u}} \cdot (\hat{w}\circ \hat{m}) ) \\ 
&  &            + (2\hat{u}\cdot (\sqrt{\hat{w}\hat{u}} \circ \hat{m})  +  
                       2\sqrt{\hat{w}\hat{u}} \cdot (\hat{u}\circ \hat{m})  ) )  =0
\end{eqnarray*}
by Lemma \ref{lemma-tr2}. $\Box$

\begin{theorem}
\label{com2}
Let  $(F, +, *)$ be a commutative presemifield.  
Then the following forms a complete set of MUBs:
$$B_{\infty}= \{ e_w \mid w \in F \}, \quad B_{m} = \{ b_{m,v} \mid v \in F \}, \quad m\in F ,$$ 
$$b_{m,v} = \frac{1}{\sqrt{q}} \sum_{w \in F} 
\omega^{{\rm Tr} (\hat{m}(\hat{w} * \hat{w}) + 2\hat{w}\hat{v})}e_w .$$
\end{theorem}

{\em Proof}. 
Let $x * y = \sum a_{ij} x^{2^i}y^{2^i}$. 
Working as in section \ref{odd-char} we see that the corresponding symplectic presemifield $(F, +, \circ )$ 
is given by multiplication:
$$x \circ y =  \sum a_{ii}^{2^{r-i}} x y^{2^{r-i}} + 
                   \sum_{i < j}a_{ij}^{2^{r-i}} x^{2^{j-i}} y^{2^{r-i}}  +  
                   \sum_{i < j}a_{ij}^{2^{r-j}} x^{2^{r+i-j}} y^{2^{r-j}}. $$
We have 
\begin{eqnarray*}
{\rm Tr} (\hat{w}\cdot (\hat{w}\circ \hat{m})) 
& = & {\rm Tr}  (  \sum \hat{w}\widehat{a_{ii}}^{2^{r-i}} \hat{w} \hat{m}^{2^{r-i}} + 
                   \sum_{i < j} \hat{w}\widehat{a_{ij}}^{2^{r-i}} \hat{w}^{2^{j-i}} \hat{m}^{2^{r-i}}  \\  
&    &       +  \sum_{i < j} \hat{w}\widehat{a_{ij}}^{2^{r-j}} \hat{w}^{2^{r+i-j}} \hat{m}^{2^{r-j}}  )  \\
& = & {\rm Tr}  (  \sum \hat{w}^{2^i}\widehat{a_{ii}} \hat{w}^{2^i} \hat{m} + 
                   \sum_{i < j} \hat{w}^{2^i}\widehat{a_{ij}} \hat{w}^{2^j} \hat{m}  +  
                   \sum_{i < j} \hat{w}^{2^j}\widehat{a_{ij}} \hat{w}^{2^i} \hat{m}  )  \\
& = & {\rm Tr}  (\hat{m}(\hat{w} *\hat{w} )).
\end{eqnarray*}

Now, the statement of theorem follows from Theorem \ref{symp2}.  $\Box$

The above theorems allow us to construct complete sets of MUBs directly from presemifields 
(in particular, from Kantor-Williams \cite{Kan5} semifields).


\subsection{Desarguesian spreads}
\label{des2}

Desarguesian spreads are constructed with the help of finite fields, and the corresponding commutative 
and symplectic semifields are the same: $x*y = x\circ y =xy$. 
Theorem \ref{autom} and \cite{Kos} imply that for the automorphism group of the 
correspon\-ding complete set of MUBs  we have 
$${\rm Aut}(\mathcal{B}) \cong F^2 . (SL_2(q) . \mathbb{Z}_r) . \mathbb{Z}_2 ,$$
where the last factor $\mathbb{Z}_2$ is absent in the case $q=2$. 
We note that the extension $F^2 . (SL_2(q) . \mathbb{Z}_r)$ is non-split 
for $q\ge 8$ (see \cite{Ab1}), but the automorphism group of the corresponding plane 
is a semidirect product.

\subsection{L\"{u}neburg planes and Suzuki groups}

In this subsection we consider non-semifield symplectic spreads related to L\"{u}neburg planes 
and Suzuki groups \cite{Kan4,Tits}. We give an explicit construction of MUBs in terms of Galois rings. 
Let $F=\Fq$ be a finite field of order $q=2^{2k+1}$ and let 
$\sigma : \alpha \rightarrow \alpha^{2^{k+1}}$ be an automorphism of  $\mathbb{F}_q$, 
$V= \mathbb{F}_q \oplus \mathbb{F}_q$. 

The symplectic spread $\Sigma$ of a space $W=V \oplus V$ consists of the subspace 
$\{(0,y) \mid y \in V\}$ and subspaces $\{(x,xM_c) \mid x \in V\}$, $c \in V$, where 

$$c=(\alpha,\beta) \in V , \ \  
M_c = \left(\begin{array}{cc}
             \alpha   &   \alpha^{\sigma^{-1}} +  \beta^{1+\sigma^{-1}}   \\
             \alpha^{\sigma^{-1}} +  \beta^{1+\sigma^{-1}}    &  \beta
      \end{array}
\right).$$

The automorphism group of this symplectic spread is 
$${\rm Aut}(\Sigma) = {\rm Aut}(Sz(q)) = Sz(q).\mathbb{Z}_{2k+1},$$
where $Sz(q) = {}^2 B_2(q)$ is the Suzuki twisted simple group.

The corresponding orthogonal decomposition of the Lie algebra $L=sl_{q^2}(\C)$ is given by 
$$L=H_{\infty} \oplus_{c \in V} H_c,$$
$$H_{\infty}= \langle D_{0,y} \mid y \in V, y\neq 0 \rangle, \quad 
H_{c} = \langle D_{x,xM_c} \mid x \in V, x\neq 0 \rangle.$$

Let $R= GR(4^{2k+1})$ be a Galois ring of characteristic 4 and cardinality $4^{2k+1}$. 
For $v=(\alpha,\beta) \in V$ we define the Teichm\"{u}ller lift of $v$ as $(\hat{\alpha},\hat{\beta})$, 
and for the matrix 
$M = \left(\begin{array}{cc}
                   \alpha      &   \beta   \\
                   \gamma   &  \delta
                \end{array}
\right)$
we define the Teichm\"{u}ller lift of $M$ as 
$\hat{M} = \left(\begin{array}{cc}
                   \hat{\alpha}      &   \hat{\beta}   \\
                   \hat{\gamma}   &  \hat{\delta}
                \end{array}
\right)$.  
Then the following bases  form a complete set of MUBs:
$$B_{\infty}= \{ e_w \mid w \in V \}, \quad B_{c} = \{ b_{c,v} \mid v \in V \}, \quad c\in V ,$$ 
$$b_{c,v} = \frac{1}{q} \sum_{w \in V} 
               \omega^{{\rm Tr} (\hat{w}\cdot \hat{w} \hat{M_c} + 2\hat{v}\cdot\hat{w})}e_w.$$
Indeed,
\begin{eqnarray*}
D_{u,uM_c}(b_{c,v}) 
& = & \frac{1}{q} \sum_{w\in V} \omega^{{\rm Tr} (\hat{w} \cdot \hat{w}\widehat{M_c} + 
                           2\hat{v}\cdot\hat{w} + 2\hat{u}\widehat{M_c}\cdot\hat{w})} e_{w+u} \\
& = & \frac{1}{q} \sum_{w\in V} \omega^{{\rm Tr} (\widehat{(w+u)} \cdot \widehat{(w+u)}\widehat{M_c} + 
                           2\hat{v}\cdot\widehat{(w+u)} + 2\hat{u}\widehat{M_c}\cdot\widehat{(w+u)})} e_{w}. 
\end{eqnarray*}
Now we set $w=(w_1,w_2)$, $u=(u_1,u_2)$. Then 
$$\widehat{(w+u)} = (\widehat{w_1} + \widehat{u_1} + 2\sqrt{\widehat{w_1}\widehat{u_1}},   
                               \widehat{w_2} + \widehat{u_2} + 2\sqrt{\widehat{w_2}\widehat{u_2}})  = 
                             \hat{w} + \hat{u} +2\hat{z}, $$
where $\hat{z} = (\sqrt{\widehat{w_1}\widehat{u_1}}, \sqrt{\widehat{w_2}\widehat{u_2}})$. 
Therefore, 
\begin{eqnarray*}
D_{u,uM_c}(b_{c,v}) 
& = & \frac{1}{q} \sum_{w\in V} 
         \omega^{{\rm Tr} ((\hat{w}+\hat{u}+2\hat{z}) \cdot (\hat{w}+\hat{u}+2\hat{z})\widehat{M_c} + 
                  2\hat{v}\cdot (\hat{w}+\hat{u}) + 2\hat{u}\widehat{M_c}\cdot (\hat{w}+\hat{u})} e_{w} \\
& = & \frac{1}{q} \sum_{w\in V} 
         \omega^{{\rm Tr} [(\hat{w} \cdot \hat{w}\widehat{M_c}+2\hat{v}\cdot \hat{w}) + 
                  (3\hat{u} \cdot \hat{u}\widehat{M_c} + 2\hat{v}\cdot \hat{u}) ]} e_{w} \\
& = & \omega^{{\rm Tr} ( 2\hat{v}\cdot \hat{u} -\hat{u}\cdot \hat{u}\widehat{M_c} )} b_{c,v},
\end{eqnarray*}
which shows that the vectors $b_{c,v}$ are common eigenvectors for the Cartan sub\-algebra 
$H_c$. 


\section{Pseudo-planar functions and MUBs}

Let $F=\Fq$ be a finite field of even order $q$.  
There are no planar functions over fields of even characteristic. 
Recently Zhou \cite{Zhou,Schm} introduced the notion of ``planar" functions in even characteristic, 
however we adopt another term and  call a function $f : F \rightarrow F$ pseudo-planar  if 
$$x \mapsto f(x+a)+f(x)+ax$$
is a permutation of $F$ for each $a\in F^*$.  
Using the Teichm\"{u}ller lift, we can also consider $f$ as a function  $f : \T \rightarrow \T$.

\begin{theorem}
\label{planar2}
Let $F=\Fq$ be a finite field of even characteristic and $f$ be a pseudo-planar function. 
Then the following forms a complete set of MUBs:
$$B_{\infty}= \{ e_w \mid w \in F \}, \quad B_{m} = \{ b_{m,v} \mid v \in F \}, \quad m\in F ,$$ 
$$b_{m,v} = \frac{1}{\sqrt{q}} \sum_{w \in F} \omega^{{\rm Tr} (\hat{m}(\hat{w}^2 + 2f(\hat{w})) + 2\hat{v}\hat{w})}e_w,$$
\end{theorem}

{\em Proof}. If $m=m_1$ then 

$$(b_{m,v}, b_{m_1,v_1}) = \frac{1}{q} \sum_{w \in F} \omega^{{\rm Tr} (2(\hat{v}-\hat{v_1})\hat{w})} =
\left\{ \begin{array}{ll} 
1 & \mbox{if $v=v'$} \\ 
0 & \mbox{if $v \neq v'$.} 
\end{array} 
\right. $$

 If $m \neq m_1$ then 
\begin{eqnarray*}
(b_{m,v}, b_{m_1,v_1}) 
& = & \frac{1}{q} \sum_{w \in F} \omega^{{\rm Tr} ((\hat{m}-\widehat{m_1})(\hat{w}^2 + 2f(\hat{w})) 
                                      + 2(\hat{v}-\widehat{v_1})\hat{w})} \\
& = & \frac{1}{q} \sum_{w \in F} \omega^{{\rm Tr} (M(\hat{w}^2 + 2f(\hat{w})) + 2\hat{u}\hat{w})}, 
\end{eqnarray*}
where $M=\hat{m}-\hat{m_1} \not\in 2R$, $\hat{u}=\widehat{v-v_1}$. Then 
\begin{eqnarray*}
|(b_{m,v}, b_{m_1,v_1})|^2 
& = & \frac{1}{q^2} \sum_{w,w_1 \in F} 
                    \omega^{{\rm Tr} (M(\hat{w}^2 + 2f(\hat{w})) + 2\hat{u}\hat{w} - 
                          M((\widehat{w_1})^2 + 2f(\widehat{w_1})) - 2\hat{u}\widehat{w_1})} \\
& = & \frac{1}{q^2} \sum_{w,w_1 \in F} \omega^{{\rm Tr} (M(\hat{w}^2 - (\widehat{w_1})^2 + 
               2f(\hat{w}) - 2f(\widehat{w_1})) + 2\hat{u}(\hat{w}-\widehat{w_1}))}.  
\end{eqnarray*}
Let $(\widehat{w_1})^2 \equiv \hat{w}^2 +\hat{a}^2 \!\!\! \pmod{2R}$, $a \in F$. 
Then $(\widehat{w_1})^2= \hat{w}^2+\hat{a}^2+2\hat{w}\hat{a}$.  Therefore, 
\begin{eqnarray*}
|(b_{m,v}, b_{m_1,v_1})|^2 
& = & \frac{1}{q^2} \sum_{w,a \in F} \omega^{{\rm Tr} (M(-\hat{a}^2-2\hat{w}\hat{a} + 
                   2f(\hat{w}) - 2f(\widehat{w+a})) - 2\hat{u}\hat{a})} \\
& = & \frac{1}{q^2} \sum_{a \in F} \omega^{{\rm Tr} (-M\hat{a}^2-2\hat{u}\hat{a})} 
                             \sum_{w \in F} \omega^{{\rm Tr} (2M(f(\widehat{w+a}) + f(\hat{w}) +\hat{w}\hat{a} ))}.  
\end{eqnarray*}
Since $f(w)$ is a pseudo-planar function, we have 
$$ \sum_{w \in F} \omega^{{\rm Tr} (2M(f(\widehat{w+a}) + f(\hat{w}) + \hat{w}\hat{a} ))} =
\left\{ \begin{array}{ll} 
q & \mbox{if $a=0$} \\ 
0 & \mbox{if $a \neq 0$.} 
\end{array} 
\right. $$
Hence,
$$|(b_{m,v}, b_{m_1,v_1})|^2 = \frac{1}{q^2} \cdot q = \frac{1}{q} . \ \ \Box$$


The following theorem shows connections between pseudo-planar functions and commutative presemifields. 

\begin{theorem}
Let $F$ be a finite field of characteristic two. 

1. If  $(F, +, *)$ is a commutative presemifield with multiplication given by 
$$x*y = xy +  \sum_{i < j}a_{ij}(x^{2^i}y^{2^j} + x^{2^j}y^{2^i})$$ 
then $f(x)=  \sum_{i < j}a_{ij}x^{2^i + 2^j}$ is a pseudo-planar function 
and $x*y = xy + f(x+y) + f(x) + f(y)$. 

2. If  $(F, +, *)$ is a commutative presemifield then there exist a strongly  isotopic commutative presemifield 
$(F, +, \star )$ and a pseudo-planar function $f$ such that $x\star y =  xy +f(x+y) + f(x) + f(y)$. 
Therefore, up to isotopism, any commutative semifield can be described by pseudo-planar 
functions. 

3. Let $f$ be a pseudo-planar function. Then  $(F, +, *)$ with multiplication $x*y = xy +f(x+y) + f(x) + f(y)$ 
is a presemifield if and only if $f$ is a quadratic function. 
\end{theorem}

{\em Proof}. 1.  If $x*y = xy +f(x+y) + f(x) + f(y)$ determines a presemifield  then the map 
$x \mapsto xy + f(x+y)+f(x)+f(y)$ is a permutation for any $y\in F^*$, therefore the map 
$x \mapsto xy + f(x+y)+f(x)$ is a permutation for any $y\in F^*$. 

2. Assume that the presemifield $(F, +, *)$ is given by 
$$x * y = \sum a_i x^{2^i}y^{2^i} + \sum_{i < j}a_{ij}(x^{2^i}y^{2^j} + x^{2^j}y^{2^i}).$$ 
Denote $g(x) = \sum a_i x^{2^i}$. Then $x*x=g(x^2)$. We need only to proof that the linear function 
$g$ is invertible and then we can take $x\star y = g^{-1}(x*y)$ and apply part 1.    

Suppose that $g$ is not invertible. Then there exists $a\in F^*$ such that $g(a^2)=0$. Then 
$$(x+a)*(x+a) = g((x+a)^2) = g(x^2 +a^2) = g(x^2) = x*x, $$
which means $x*x + x*a + a*x + a*a = x*x$ and $a*a =0$, a contradiction. 

2. The condition $(x+z)*y = x*y + z*y$ is equivalent to 
$$(x+z)y+f(x+z+y) + f(x+z) + f(y) $$
$$= xy+f(x+y) + f(x) + f(y) + zy+f(z+y) + f(z) + f(y).$$
Therefore, 
$$f(x+y+z)+f(x+y)+f(x+z)+f(y+z)+f(x)+f(y)+f(z)=0, $$
which means that $f$ is a quadratic function. $\Box$

{\bf Remark 2}. 
The notion of pseudo-planar functions can be defined over a field $F$ of arbitrary characteristic. 
We call a function $f : F \rightarrow F$ pseudo-planar  if the map 
$x \mapsto f(x+a)-f(x)+ax$
is a permutation of $F$ for each $a\in F^*$.  Such functions carry similar properties 
as pseudo-planar functions in even characteristic (inclu\-ding MUBs constructions). 
In particular, 
if $f$ is a quadratic pseudo-planar function then the product $x*y = xy + f(x+y)-f(x)-f(y)$ 
defines a commutative presemifield $(F,+,*)$. 
Probably such a definition of pseudo-planar functions provides unified approach to all characteristics, 
at least it provides a bridge between the even and odd characteristic cases. 
Note that in the case of odd characteristic function $f$ is pseudo-planar if and only if 
the function $x^2+2f(x)$ is planar, and formulas in Theorems \ref{planar} and \ref{planar2} 
have a unified look in the language of pseudo-planar functions.


\bigskip

{\bf Acknowledgments} 

\medskip

The author would like to thank William M. Kantor, Claude Carlet and Yue Zhou for valuable discussions 
on the content of this paper. 
This research was supported by UAEU grant 21S073 and NRF grant 31S088. 



\begin{thebibliography}{99}



\bibitem{Ab1} 
K. S. Abdukhalikov, 
On a group of automorphisms of an orthogonal decomposition of a Lie algebra of type $A_{2^m-1}$,  
(Russian) Izv. Vyssh. Uchebn. Zaved. Mat., no. 10, 11--14 (1991); 
English translation in Soviet Math. (Iz. VUZ) 35,  9--12 (1991).

\bibitem{Ab2} 
K. S. Abdukhalikov, 
Invariant integral lattices in Lie algebras of type $A_{p^m-1}$, 
(Russian) Mat. Sb. 184, no. 4, 61--104 (1993); 
English translation in Russian Acad. Sci. Sb. Math. 78, no. 2, 447--478 (1994). 

\bibitem{Ab3} 
K. S. Abdukhalikov, 
Affine invariant and cyclic codes over $p$-adic numbers and finite rings, 
Des. Codes Cryptogr. 23,  343--370 (2001).

\bibitem{Ab4}
K. Abdukhalikov, E. Bannai and S. Suda, 
Association schemes related to universally optimal configurations, 
Kerdock codes and extremal Euclidean line-sets, 
J. Combin. Theory Ser. A 116, 434--448 (2009).

\bibitem{All}
W. O. Alltop, 
Complex sequences with low periodic correlations, 
IEEE Trans. Inform. Theory 26 (3), 350--354 (1980).

\bibitem{Bad}
L. Bader, W. M. Kantor and G. Lunardon,
Symplectic spreads from twisted fields, 
Boll. Un. Math. Ital. A (7) 8,  383--389 (1994). 

\bibitem{Ball}
S. Ball, J. Bamberg, M. Lavrauw and T. Penttila, 
Symplectic spreads, 
Des. Codes Cryptogr. 32, 9--14 (2004).

\bibitem{Ban}
S. Bandyopadhyay, P. O. Boykin, V. Roychowdhury and F. Vatan, 
A new proof for the existence of mutually unbiased bases, 
Algorithmica 34,  512--528 (2002).

\bibitem{Ben}
C. H. Benneth and G. Brassard, 
Quantum Cryptography: Public key distribution and coin tossing, 
in Proceedings of the IEEE International Conference on Computers, Systems, and Signal Processing, 
Bangalore  175 (1984). 

\bibitem{Bier}
J. Bierbrauer, 
New semifields, PN and APN functions,  
Des. Codes Cryptogr. 54(3), 189--200 (2010).

\bibitem{Boy}
P. O. Boykin, M. Sitharam, Pham Huu Tiep and P. Wocjan,
Mutually unbiased bases and orthogonal decompositions of Lie algebras,
Quantum Inf. Comput. 7, no. 4, 371--382 (2007).

\bibitem{Buda}
L. Budaghyan and T. Helleseth, 
New perfect nonlinear multinomials over $F_{p^{2k}}$ for any odd prime $p$. 
In: SETA ’08: Proceedings of the 5th International Conference on Sequences and their
Applications, 403--414. Springer, Berlin, Heidelberg (2008).

\bibitem{Cal}
A. R. Calderbank, P. J. Cameron, W. M. Kantor and J. J. Seidel,
$\mathbb{Z}_4$-Kerdock codes, orthogonal spreads, and extremal Euclidean line-sets,
Proc. London Math. Soc. (3) 75, 436--480 (1997).

\bibitem{Coh}
S. D. Cohen, M. J. Ganley, 
Commutative semifields, two-dimensional over their middle nuclei, 
J. Algebra 75, 373--385 (1982).

\bibitem{Coul}
R. S. Coulter and R. W. Matthews, 
Planar functions and planes of Lenz-Barlotti class II,  
Des. Codes Cryptogr. 10, 167--184 (1997).

\bibitem{Coul2}
R. S. Coulter, M. Henderson and P. Kosick, 
Planar polynomials for commutative semifields with specified nuclei, 
Des. Codes Cryptogr. 44,  275--286 (2007).

\bibitem{Coul3}
R. S. Coulter and P. Kosick, 
Commutative semifields of order 243 and 3125, 
in: Finite Fields: Theory and Applications,
in: Contemp. Math., vol. 518, Amer. Math. Soc., Providence, RI, 129--136 (2010).

\bibitem{Dem}
P. Dembowski, 
Finite geometries, Springer, Berlin (1968). 

\bibitem{Dick}
L. E. Dickson, 
On commutative linear algebras in which division is always uniquely possible,  
Trans. Amer. Math. Soc. 7, 514--522 (1906).

\bibitem{Ding1}
C. Ding and J. Yuan, 
A family of skew Hadamard difference sets,  
J. Combin. Theory Ser. A 113, 1526--1535 (2006).

\bibitem{Ding}
C. Ding and J. Yin, 
Signal sets from functions with optimum nonlinearity, 
IEEE Trans. Communications 55, No. 5,  936--940 (2007).

\bibitem{Gan}
M. J. Ganley, 
Central weak nucleus semifields, 
European J. Combin. 2, 339--347 (1981).

\bibitem{God}
C. Godsil and A. Roy, 
Equiangular lines, mutually unbiased bases, and spin models, 
European J. Combin. 30, 246--262 (2009). 

\bibitem{Gow} 
R. Gow, 
Generation of mutually unbiased bases as powers of a unitary matrix in 2-power dimensions, 
arXiv:math/0703333.

\bibitem{Ham}
R. Hammons, P. V. Kumar, A. R. Calderbank, N. J. A. Sloane and P. Sol\'{e},
The $Z_4$-linearity of Kerdock, Preparata, Goethals, and related codes,
IEEE Trans. Inform. Theory 40, No. 2, 301--319 (1994).

\bibitem{Iva}
I. D. Ivanovi\'{c}, 
Geometrical description of quantal state determination,  
J. Phys. A 14,  3241--3245 (1981).

\bibitem{Kan}
W. M. Kantor and M. E. Williams, 
New flag-transitive affine planes of even order, 
J.~Comb. Theory(A) 74,  1--13 (1996).


\bibitem{Kan8}
W. M. Kantor,
Ovoids and translation planes, 
Canad. J. Math. 34, 1195--1207 (1982).

\bibitem{Kan9}
W. M. Kantor,
Projective planes of order q whose collineation groups have order $q^2$,  
J.~Alg. Comb. 3, 405--425 (1994).


\bibitem{Kan6}
W. M. Kantor,
Note on Lie algebras, finite groups and finite geometries, 
in: Groups, Difference Sets, and the Monster (Eds. K. T. Arasu et al.)  
73--81, de Gruyter, Berlin-New York (1996).

\bibitem{Kan3}
W. M. Kantor,
Commutative semifields and symplectic spreads, 
J. Algebra 270, 96--114 (2003).

\bibitem{Kan5}
W. M. Kantor and M. E. Williams, 
Symplectic semifield planes and $\Zf$-linear codes,  
Trans. Amer. Math. Soc. 356, 895--938 (2004).

\bibitem{Kan2}
W. M. Kantor and M. E. Williams, 
Nearly flag-transitive affine planes,  
Adv. Geom. 10, 161--183 (2010).

\bibitem{Kan4}
W. M. Kantor,
MUBs inequivalence and affine planes, 
J. Math. Phys. 53, 032204 (2012).

\bibitem{Klap} 
A. Klappenecker and M. R\"{o}tteler, 
Constructions of mutually unbiased bases, 
in: Finite Fields and Applications, Lecture Notes in Computer Sciences, 
vol. 2948, 137--144, Springer, Berlin (2004). 

\bibitem{Knu}
D. E. Knuth, 
Finite semifields and projective planes, 
J. Algebra 2, 182--217 (1965). 

\bibitem{Kos}
A. I. Kostrikin and Pham Huu Tiep,
Orthogonal decompositions and integral lattices,
Walter de Gruyter, Berlin (1994).

\bibitem{LeC}
N. LeCompte, W. J. Martin and W. Owens, 
On the equivalence between real mutually unbiased bases and a certain class of association schemes, 
European J. Combin., 31, 1499--1512 (2010).

\bibitem{Nech}
A. A. Nechaev, 
Kerdock code in a cyclic form, 
(Russian) Diskret. Mat. 1, no. 4, 123--139 (1989); 
English translation in Discrete Math. Appl. 1, no. 4, 365--384 (1991).  

\bibitem{Pen}
T. Penttila and B. Williams, 
Ovoids of parabolic spaces. Geom. Dedicata 82, 1--19 (2000).

\bibitem{Roy}
A. Roy and A. J. Scott, 
Weighted complex projective 2-designs from bases: optimal
state determination by orthogonal measurements, 
J. Math. Phys. 48, 072110  (2007).

\bibitem{Schm}
K.-U. Schmidt and Y. Zhou, 
Planar functions over fields of characteristic two, 
J. Algebraic Combin. 40, 503--526 (2014). 
 
\bibitem{Schw}
J. Schwinger, Unitary operator bases, 
Proc. Nat. Acad. Sci. U.S.A. 46, 570-579 (1960).

\bibitem{Sue}
C. Suetake, A new class of translation planes of order $q^3$, 
Osaka J. Math. 22, 773--786 (1985).

\bibitem{Thas}
J. A. Thas and S. E. Payne, 
Spreads and ovoids in finite generalized quadrangles,  
Geom. Dedicata 52, 227--253 (1994).

\bibitem{Tits} 
J. Tits, Ovo\"{\i}des et groupes de Suzuki,  
Arch. Math. 13, 187--198 (1962).

\bibitem{Weng}
G. Weng and X. Zeng, 
Further results on planar DO functions and commutative semifields, 
Des. Codes Cryptogr. 63, 413--423 (2012).

\bibitem{Woo}
W. K. Wootters and B. D. Fields, 
Optimal state-determination by mutually unbiased measurements, 
Ann. Physics 191 (2), 363--381 (1989). 

\bibitem{Zha}
Z. Zha, G. M. Kyureghyan and X. Wang, 
Perfect nonlinear binomials and their semifields, 
Finite Fields and Their Applications 15, 125--133 (2009).

\bibitem{Zha2}
Z. Zha and X. Wang, 
New families of perfect nonlinear polynomial functions, 
J. Algebra 322, 3912--3918 (2009). 

\bibitem{Zhou}
Y. Zhou, 
$(2^n,2^n,2^n,1)$-relative difference sets and their representations,  
J. Combin. Designs 21 (12), 563--584 (2013).

\bibitem{Zhou2}
Y. Zhou and A. Pott, 
A new family of semifields with 2 parameters, 
Adv. Math. 234, 43--60 (2013).

\end{thebibliography}
\end{document}